# Contribution of Non Integer Integro-Differential Operators (NIDO) to the geometrical understanding of Riemann's conjecture-(II)

Alain Le Méhauté*
Abdelaziz El Kaabouchi
Laurent Nivanen
Institut Supérieur des Matériaux et Mécaniques Avancés (ISMANS)
Campus au Mans de l'Université du Québec
72 000 Le MANS
France
alm@ismans.fr

*Abstract* – Advances in fractional analysis suggest a new way for the physics understanding of Riemann's conjecture. It asserts that, if s is a complex number, the non trivial zeros of zeta function $\frac{1}{\zeta(s)} = \sum_{n=1}^{\infty} \frac{\mu(n)}{n^s}$ in the gap [0,1], is characterized by $s = \frac{1}{2}(1+2i\theta)$. This conjecture can be understood as a consequence of 1/2-order fractional differential characteristics of automorph dynamics upon opened punctuated torus with an angle at infinity equal to $\pi/4$. This physical interpretation suggests new opportunities for revisiting the cryptographic methodologies.

*Keywords* - Fractal; Algebraic structures and number theory; Differential geometry and topology; Statistical mechanics; Fractal analysis; Signal treatment; Cryptography.

*PACS* : O5.45.df, 02.10.De, 02.40.-k,05.20.-y

In the huge amount of scientific works left by Euler in 1737 we find the following relation $\sum_{n=1}^{\infty} \frac{1}{n^s} = \prod_{p} \left(1 - \frac{1}{p^s}\right)^{-1}$, for real $s > 1$, where the product is taken over all primes. In 1859 Riemann published a short revolutionary note [1,2] introducing the function $\zeta(s)$ by $\frac{1}{\zeta(s)} = \sum_{n=1}^{\infty} \frac{\mu(n)}{n^s}$, where $s > 0$ ($\mu$ is the Möbius function ($\mu(1) = 1$ and $\mu(n) = (-1)^k$ if $n$ is the product of $k$ different prime numbers, and $\mu(n) = 0$ in other cases). He assumed that non trivial zeros of $\zeta(s)$ are located on the straight line given by $s = \sigma[1 \pm (1/\sigma)i\theta]$ with $\sigma = 1/2$. The physical meaning of this conjecture, as well as the mathematical arguments supporting the general assumption [3,4,5,6,7] remains unknown, in spite of the Siegel's historical studies [3], and the Weil's resolution of Riemann's conjecture in finite field [4].

In this note, we address a constructivist proposal aiming at the physical understanding of Riemann's conjecture. This proposal is based on the analytical parameterization of transfer function of non integer Integro-Diffential Equations [8,9,10] in the set of natural numbers when the tilling of the complex plane is reduced in the set of prime numbers [11].

The origin of the present work is an analogy between the pinning of the $\sigma$ value on the fraction 1/2 and the pinning of phase angle upon the value $\varphi = \pm \pi/4$ for any dynamics which is under the control of 1/2-order Differential Equation [8,9,10]. Our analysis will be carried out in two major steps:

*1) First assertion:* we can built two classes of $\eta(s) = (1 - 2^{1-s})\zeta(s)$ function with $0 < s < 1$ from the parameterization of the class of Cole and Cole transfer functions [8,9,12] with $\sigma = 1/d$ and $s(d,\theta) = \sigma[1 + (1/\sigma)i\theta]$ *within the set of natural numbers* **N**. The geometrical meaning of this function is the sum of hyperbolic distances upon the automorph geodesics of the opened punctuated torus with an angle at infinity equal to $\varphi = (\pi/2)\sigma$. This class of functions can only be equal all in the complex plane over both modulus and phases if and only if $\sigma = 1/2$ [11]. This property might explain the position of non trivial zeros for $\zeta(s)$.

*2) Second assertion:* when the problem is moved from **N** to **P** (set of all prime numbers), analytical prolongation of $\zeta(s)$ function is found based on the parametrization of transfer functions.

---

Now let us consider in detail both assertions. First let us recall some well established results about fractional dynamics associated with geodesics on opened punctuated torus [8,9,11]. We call Cole and Cole transfer function [12] a function Z(v):

$$Z(v) = \frac{Z_0}{1 + \left[\frac{iv}{v_c}\right]^{1/d}} \qquad (1)$$

This arc in the complex plane, $1/d$ hyperbolic geodesic $H_\varphi(p_i)$, may be analytically written as a part of a circle which centre is located out of the abscise axis [8,9]; This function is itself associated with a very simple fractional differential equation based on any intensive and extensive thermodynamic variables : U(t), I(t)

$$(1/v_c)^{1/d} \frac{d^{1/d}}{dt^{1/d}}[U(t)] = Z_0 I(t) - U(t) \qquad (2)$$

where $v$ is the Fourier transform of a time variable $t$ which parameterises the dynamics, and $\varphi$ is the pinned phase angle given by $\varphi = \frac{\pi}{2}(1 - \frac{1}{d})$ leading the azimuth of the centre of the support circle of the arc. We call $p_i$ the length of the string (abscise) of the arc $H_\varphi(p_i)$. The main feature of $Z(v)$ (Fig.1) is the following: whatever $p_i$ the Z-parameterization requires only the class of hyperbolic distances $\delta$ parameterised by a variable $v_c$ and $d$ [8,9,10] such as $\delta^\rightarrow = \frac{u}{v} = \left(\frac{v}{v_c}\right)^{1/d}$. As shown in Fig.1, in the complex plan the definition of $p_i^\rightarrow$ involves the definition of another string $p_j^\downarrow$ on the imaginary axis (Fig.1), both arcs being analytically related by the same class of hyperbolic distance $\delta^\downarrow = \frac{u'}{v'} = \left(\frac{v}{v_c}\right)^{1/D}$ with $\frac{1}{D} = 1 - \frac{1}{d}$. More precisely, if we reduced the set of parameters to the natural number **N** instead of the rational **Q** and using symmetric properties of hyperbolic distances, these distances may split into four classes with indices $l$ (left) and $r$ (right) according to the position with respect to the characteristic time and position, horizontal $(\rightarrow)$ or vertical $(\downarrow)$ in such a way that

$$\delta_l^\rightarrow = \frac{u}{v} = \left(\frac{v}{v_c}\right)^{1/d} = (n)_l^{1/d}, \quad \delta_r^\rightarrow = \frac{u}{v} = \left(\frac{v}{v_c}\right)^{1/d} = (1/n)_r^{1/d},$$

$$\delta_r^\downarrow = \frac{u'}{v'} = \left(\frac{v'}{v'_c}\right)^{1/D} = (n)_r^{1/D}, \quad \delta_l^\downarrow = \frac{u'}{v'} = \left(\frac{v'}{v'_c}\right)^{1/D} = (1/n)_l^{1/D}$$

From that starting point we can adopt a constructivist approach described in chart 1 where a function $\eta$ is constructed step by step via the definition of hyperbolic distance (step1), the parameterization of this distance using a multiplication with a complex function of free $\theta$ phase angle (step 2), and finally via the building of two 'finite summated distances' $\lambda$ in the natural number set (step 3&4).
We put $s_1 = (1/d)(1 + di\theta)$ and $s_2 = (1/D)(1 + Di\theta)$.

Both functions $\eta$ are related to zeta Riemann function according to $\zeta(s)^\rightarrow = \zeta(s)^\downarrow = (1 - 2^{1-s})\eta^\downarrow(s) = (1 - 2^{1-s})\eta^\rightarrow(s)$ [3,4]. Zero of $\zeta$ function gives also zeros of $\eta$ functions. If both zeros are different there is a relation between $d$ and $\theta$ such as $\frac{1}{d(\theta)} = 1 - \frac{1}{D(\theta)}$ but also in unique complex plane $\frac{1}{D(\theta)} = \frac{1}{d(\theta)}$ therefore $d$ is a constant and $\frac{1}{d(\theta)} = \frac{1}{2}$. All zero of zeta Riemann function are located on the vertical line with an abscise $1/2$ and $s = (1/2)(1 + 2i\theta)$.

Above constraints involves $p_i = p_j$ [11]. This relation leads to the second assertion, the relation between $\zeta(s) = 0$ and the set of prime number **P**. To understand this relation let us come back to the foundation of the dynamic on fractal geometry. Briefly speaking, as shown by Mandelbrot 30 years ago [14], the metric of fractal in 2D is characterized by the Mandelbrot's equation

$(1/\delta)^d p = v_c$ where $1/\delta$ is the gauge of geometrical measure, $d \in \mathbf{R}$ is the non integer dimension, $p$ is a detailed account of the number of gauges (balls in 2D) required to cover the fractal set, and $v_c$ is a constant -the d-Hausdorff's content-, independent of $\delta$. The physical 'component' matching any dissipative fractional dynamics is a *bounded folded curvature capacitor* initially called '*Fractance*' [8] for 'capacitor on fractal geometry'. Due to space-time coupling, a screening of $p$, which mimics the Fourier (Laplace) analysis, leads to a screening along the characteristic lengths $1/\delta$ (practically the Fourier Transform of the velocity with respect to the time). According to the non integer fractal metric we can write the so called Mandelbrot's relation [14]: $1/\delta - \left(\dfrac{v_c}{p}\right)^{1/d} = 0$.

The relation between the second assertion and the first one comes from the fact that for all $\delta$, $v_c$ and $p$, an exchange of energy upon *fractance* leads to a pinning of phase angle $\varphi$ according to $\varphi = \left(1 - \dfrac{1}{d}\right)\dfrac{\pi}{2}$ (Fig.1). Obviously, the screening via $p$, in the set of natural numbers $\mathbf{N}$, meets on the way all prime numbers $\mathbf{P}$. If the set of $p$ is reduced to the set of prime number $\mathbf{P}$ that is $p = p_i \in \mathbf{P}$, -that can be done without any loss of information upon the class of geodesics- assertion 1 assures that the constraints of parameterisation and tilling must be restricted to $\psi = \pi/4$ and $d = 2$, therefore $1/\delta - \left(\dfrac{v_c}{p}\right)^{1/2} = 0$.

*Chart 1: constructive approach of zeta Riemann function from hyperbolic distances*

| $1/\delta_l^{\rightarrow} = (n)^{-(1/d)}$ | $\delta_r^{\rightarrow} = (1/n)^{(1/d)}$ | $1/\delta_r^{\downarrow} = (n)^{-(1/D)}$ | $\delta_l^{\rightarrow} = (1/n)^{(1/D)}$ | Multiplication |
|---|---|---|---|---|
| $n^{i\theta}$ | $n^{-i\theta}$ | $n^{i\theta}$ | $n^{-i\theta}$ | Parametrization |
| $1/\Delta_g^{\rightarrow} = (n)^{-s_1}$ | $\Delta_r^{\rightarrow} = (n)^{s_1}$ | $1/\Delta_r^{\downarrow} = (n)^{-s_2}$ | $\Delta_l^{\downarrow} = (n)^{s_2}$ | Summation |
| $1/\xi_0(s_1)^{\rightarrow} = \sum_{n=1}^{\infty}(n)^{-s_1}$ | $\lambda_0(s_1)^{\rightarrow} = \sum_{n=1}^{\infty}(n)^{s_1}$ | $1/\xi_0(s_2)^{\downarrow} = \sum_{n=1}^{\infty}(n)^{-s_2}$ | $\lambda_0(s_2)^{\downarrow} = \sum_{n=1}^{\infty}(n)^{s_2}$ | $s = \sigma[1 \pm (1/\sigma)i\theta)]$ |
| | $\eta(s_1) = \sum_{n=1}^{\infty}\dfrac{(-1)^{n+1}}{(n)^{s_1}}$ | | $\eta(s_2) = \sum_{n=1}^{\infty}\dfrac{(-1)^{n+1}}{(n)^{s_2}}$ | Convergence for $0 < \text{Re}(s_k) < 1$ |

According to the geometrical meaning of the dynamic parameters, and due to the non divisibility of the prime number $p_i$ except by the unit, we must distinguish and identify two hypotheses concerning the Hausdorff content, i.e., either (i) the trivial solution of the previous equation given by the norm upon $1/\delta = 1$ and $v_c = p_i$ or (ii) according to the fact that in the complex plane $1^2 = \left(p_i^{i2\theta'}\right)\left(p_i^{-i2\theta'}\right)$ the introduction of a free parameter $\theta'$ to defined 'Hausdorff normalisation' $v_c = 1^2$. The last identities lead to a complex $\theta'$–parametrization of the Mandelbrot relation:

$$\dfrac{1/\delta}{p_i^{2i\theta'}} - \left(\dfrac{1}{p_i}\right)^{(1+2i\theta')/2} = 0 \text{ and } \dfrac{1/\delta}{p_i^{-2i\theta'}} - \left(\dfrac{1}{p_i}\right)^{(1-2i\theta')/2} = 0.$$

According to references [8,9,12] extended to complex dimension [13], these equations disclose a clear understanding of the role of the complex extension of real *s* via $\theta'$: when the complex plan is tilled by unit square and (equivalently) isomorph prime number square, $\pm\theta'$ must be chosen in such a way that, due to the symmetry, whatever $p_i \in \mathbf{P}$, the 'real' part $(\rightarrow)$ of the first term of above equation, $\left(1/\delta p_i^{2i\theta'}\right) + \left(1/\delta p_i^{-2i\theta'}\right)$, must be equal to unity whereas the symmetric 'imaginary' parts $(\downarrow)$ becomes nil. The parameters $\pm\theta'$ control the location of all non trivial zeros of above function. Its choice closes the simplest physical analysis of the zeros of a certain function $\varpi(s)$ which can be written under the following mathematical form

$$\varpi(s) = \prod_i \left(1 - \left[\dfrac{1}{p_i}\right]^{s(+\theta)} + \left[\dfrac{1}{p_i}\right]^{s(-\theta)}\right)$$

with $s(\theta') = (1 \pm 2i\theta')/2$.

CONCLUSION

From the physics of complex media, the present fractional analysis of Cole and Cole transfer function introduced a non trivial point of view about Riemann's conjecture: *the real part of 's' position of the zeros and the analytical expansion of ζ(s) in the gap [0,1]*. The factor ½ would be related to the singular status of ½ transfer function in the complex plane. The assumption of the existence of a relation between discrete parameterisation in $\mathbf{N}^2$ and $\mathbf{P}^2$ of 1/2- hyperbolic geodesic, and tilling of complex plan with fundamental domains associated with $\pi/4$ punctuated torus might seems conjectural but it clearly introduces at least an interesting opportunities for a constructive approach of Riemann zeta function.

It is clear that the above physical interpretation is built on the cross road of very academic mathematical topics such as hyperbolic space, Riemann manifolds, symbolic dynamics, auto-morph functions, holomorph analysis, modular algebra, singularity and measure theory, subjects which can now be related to above analysis.

Due to the prominence of the zeta function in signal treatment, the fractional point of view might shed new lights upon NP problem considered especially for applications in cryptography, coding, signal treatment and image analysis.

ACKNOWLEDGMENTS

The authors would like to thank Claude Tricot, Michel Mendes France, Jean Vaillancourt, Serge Perrine Françoise Lippi, and R.R. Nigmatullin of the Kazan State University (Russia) for partnership in theoretical physics as well as Materials Design sa. & inc and ADES Conseil for financial support.

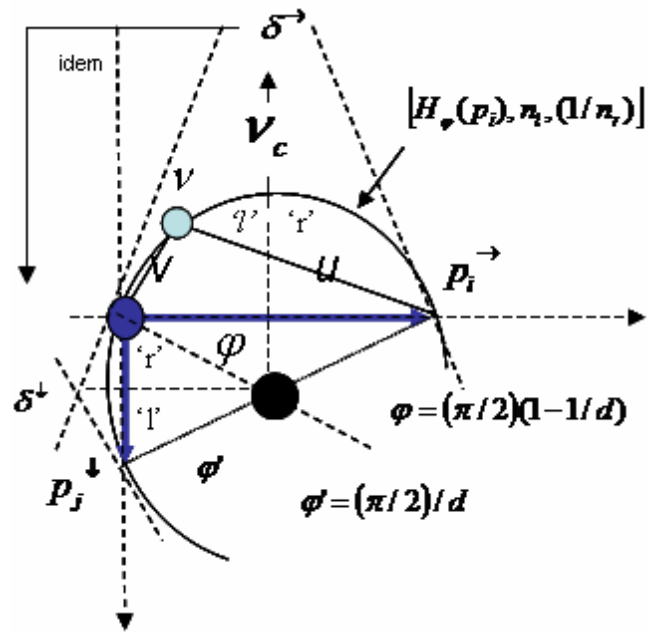

Fig 1: Representation of the different parameters given in the text analysed upon hyperbolic Cole and Cole geodesics in the complex plane. Zeros of zeta Riemann function are obtained when $\varphi = \pi/4$ and $p_i = p_j$.